\documentclass[12pt]{article}
\usepackage[cp1251]{inputenc}
\usepackage[russian]{babel}
\usepackage{amssymb}
\usepackage{amsmath}

\begin{document}

UDC 517.51


\quad

\begin{center}
\textbf{\large Uniformly convergent Fourier series\\ and
multiplication of functions}
\end{center}

\begin{center}
\textsc{V. V. Lebedev}
\end{center}

\begin{quotation}
{\small Let $U(\mathbb T)$ be the space of all continuous
functions on the circle $\mathbb T$ whose Fourier series
converges uniformly. Salem's well-known example shows that a
product of two functions in $U(\mathbb T)$ does not always
belongs to $U(\mathbb T)$ even if one of the factors belongs to
the Wiener algebra $A(\mathbb T)$. In this paper we consider
pointwise multipliers of the space $U(\mathbb T)$, i.e., the
functions $m$ such that $mf\in U(\mathbb T)$ whenever $f\in
U(\mathbb T)$. We present certain sufficient conditions for a
function to be a multiplier and also obtain some results of
Salem type.

}
\end{quotation}

\quad

\begin{center}
\textbf{1. Introduction}
\end{center}

We consider functions $f$ on the circle $\mathbb T=\mathbb
R/(2\pi\mathbb Z)$ and their Fourier expansions
$$
f(t)\sim\sum_{k\in\mathbb Z} \widehat{f}(k)e^{ikt}.
\eqno(1)
$$
Here $\mathbb R$ is the real line and $\mathbb Z$ is the
additive group of integers (we naturally identify the functions
on $\mathbb T$ with the $2\pi$-periodic functions on $\mathbb
R$).

Let $C(\mathbb T)$ be the space of all continuous functions $f$
on $\mathbb T$ with the usual norm $\|f\|_{C(\mathbb
T)}=\sup_{t\in\mathbb T}|f(t)|$. Consider the space $U(\mathbb
T)$ of all functions $f\in C(\mathbb T)$ which have  uniformly
convergent Fourier series, i.e., satisfy the condition
$\|f-S_N(f)\|_{C(\mathbb T)}\rightarrow 0$ as
$N\rightarrow\infty$, where $S_N(f)$ is the $N$th partial sum
of the Fourier series of $f$:
$$
S_N(f)(t)=\sum_{|k|\leq N} \widehat{f}(k)e^{ikt}.
$$
The norm on $U(\mathbb T)$ is defined by
$$
\|f\|_{U(\mathbb T)}=\sup_N\|S_N(f)\|_{C(\mathbb T)}.
$$
We note that the space $U(\mathbb T)$ is a Banach space.
Clearly, the inclusion $U(\mathbb T)\subseteq C(\mathbb T)$
holds with the corresponding relation for the norms:\!
$\|\cdot\|_{C(\mathbb T)}\leq\|\cdot\|_{U(\mathbb T)}$.

Another well-known space related to expansion (1) is the space
$A(\mathbb T)$ of all functions $f\in C(\mathbb T)$ whose
Fourier series converges absolutely. Endowed with the natural
norm
$$
\|f\|_{A(\mathbb T)}=\sum_{k\in\mathbb Z}|\widehat{f}(k)|,
$$
the space $A(\mathbb T)$ is a Banach space. Moreover, under the
usual multiplication of functions $A(\mathbb T)$ is a Banach
algebra (known as the Wiener algebra). Obviously, $A(\mathbb
T)\subseteq U(\mathbb T)$ and $\|\cdot\|_{U(\mathbb
T)}\leq\|\cdot\|_{A(\mathbb T)}$.

It is known that the space $U(\mathbb T)$, unlike $A(\mathbb
T)$ and $C(\mathbb T)$, is not an algebra; the product of two
functions in $U(\mathbb T)$ does not necessarily belong to
$U(\mathbb T)$. Moreover, there exist $f\in A(\mathbb T)$ and
$g\in U(\mathbb T)$ such that $fg\notin U(\mathbb T)$. This
result is due to Salem; see [1, Ch.~1, Sec.~6]. In this paper
we present certain sufficient conditions for a function $m$ to
have the property that $mf\in U(\mathbb T)$ whenever $f\in
U(\mathbb T)$. We also obtain some results of Salem type.

\quad

\begin{center}
\textbf{2. Pointwise multipliers: Sufficient conditions}
\end{center}

Let $m\in C(\mathbb T)$. We say that $m$ is a pointwise
multiplier of the space $U(\mathbb T)$ if for every function
$f$ in $U(\mathbb T)$ the product $mf$ is in $U(\mathbb T)$ as
well. We denote the space of all these multipliers by
$MU(\mathbb T)$. Clearly, if $m\in MU(\mathbb T)$, then the
operator $f\rightarrow mf$ of multiplication by $m$ is a
bounded operator on $U(\mathbb T)$. The space $MU(\mathbb T)$
endowed with the natural norm
$$
\|m\|_{MU(\mathbb T)}=\sup_{\|f\|_{U(\mathbb T)}\leq 1}\|mf\|_{U(\mathbb T)}
$$
 is a Banach algebra.
Obviously, $MU(\mathbb T)\subseteq U(\mathbb T)$.

Recall a classical fact (see, e.g., [2, Ch.~VIII, Sec.~1,
Theorem~3]): if a function $g\in C(\mathbb T)$ satisfies the
uniform Dini condition
$$
\sup_{t}\int_{|t-\theta|\leq
\varepsilon}\frac{|g(t)-g(\theta)|}{|t-\theta|}d\theta=o(1),
\qquad\varepsilon\rightarrow+0,
\eqno(2)
$$
then $g\in U(\mathbb T)$.

We shall prove the following theorem.

\quad

\textsc{Theorem 1.} \emph{If a function $m\in C(\mathbb T)$
satisfies the uniform Dini condition, then $m\in MU(\mathbb
T)$. Furthermore, the following estimate holds
$$
\|m\|_{MU(\mathbb T)}\leq\|m\|_{C(\mathbb T)}+
\sup_{t}\int_{|t-\theta|\leq\pi}\frac{|m(t)-m(\theta)|}{|t-\theta|}d\theta.
$$}

\quad

Theorem 1 immediately implies a sufficient condition for a
function to belong to $MU(\mathbb T)$ in terms of the modulus
of continuity; namely, we obtain the following theorem.

\quad

\textsc{Theorem 2.} \emph{If $m\in C(\mathbb T)$ and
$$
\int_0^\pi\frac{\omega(m, \delta)}{\delta}d\delta<\infty,
$$
where $\omega(m, \delta)=\sup_{|t_1-t_2|\leq
\delta}|m(t_1)-m(t_2)|$, then $m\in MU(\mathbb T)$.
Furthermore, one has
$$
\|m\|_{MU(\mathbb T)}\leq\|m\|_{C(\mathbb T)}+
2\int_0^\pi\frac{\omega(m, \delta)}{\delta}d\delta.
$$}

\quad

The following theorem, which gives a sufficient condition in
terms of the Fourier transform (proved by the author as a lemma
in [3]) is an immediate consequence of Theorem 2.

\quad

\textsc{Theorem 3.} \emph{If $m\in C(\mathbb T)$ and
$$
\sum_{k\in\mathbb Z}|\widehat{m}(k)|\log(|k|+2)<\infty,
$$
then $m\in MU(\mathbb T)$. Furthermore,
$$
\|m\|_{MU}\leq c \sum_{k\in\mathbb
Z}|\widehat{m}(k)|\log(|k|+2),
$$
where $c>0$ does not depend on $m$.}

\quad

To derive Theorem 3 from Theorem 2, it  suffices to note that
$$
\omega(m, \delta)\leq \sum_{k\neq 0}|\widehat{m}(k)||2\sin (k\delta/2)|
$$
and use the obvious relation\footnote{We write $a(n)\lesssim
b(n)$ or $b(n)\gtrsim a(n)$ in the case when $a(n)\leq c b(n)$
for all sufficiently large $n$, where $c>0$ does not depend on
$n$. The relation $a(n)\simeq b(n)$ means that $a(n)\lesssim
b(n)$ and $b(n)\lesssim a(n)$ simultaneously. }
$$
\int_0^\pi\frac{|\sin (k\delta/2)|}{\delta}d\delta\simeq\log |k|.
$$

\quad

\textsc{Proof of Theorem 1.} For each $N=0, 1, 2, \ldots$ let
$Q_N$ be the commutator of the operator of multiplication by
$m$ and the partial sum operator $S_N$, i.e.,
$$
Q_N : f\rightarrow mS_N(f)-S_N(mf).
$$
Considering these commutators as operators on $C(\mathbb T)$,
let us show that (under the assumptions of the theorem on $m$)
the following two conditions hold:

\quad

(i) the sequence of the norms $\|Q_N\|_{C(\mathbb T)\rightarrow
C(\mathbb T)}, \,N=0, 1, 2, \ldots,$ is bounded;

(ii) for each $n\in\mathbb Z$, the sequence of the norms
$\|Q_Ne_n\|_{C(\mathbb T)}, \,N=0, 1, 2, \ldots,$ of the images
of the exponential function $e_n(t)=e^{int}$ tends to zero as
$N\rightarrow\infty$.

\quad

We set
$$
c=\sup_{t}\int_{|t-\theta|\leq\pi}\frac{|m(t)-m(\theta)|}{|t-\theta|}d\theta.
$$
For $0<|x|\leq\pi$, the Dirichlet kernel $D_N(x)=\sum_{|k|\leq
N}e^{ikx}$ is estimated as
$$
|D_N(x)|=\bigg|\frac{\sin(N+1/2)x}{\sin(x/2)}\bigg|\leq\frac{1}{|\sin(x/2)|}\leq\frac{\pi}{|x|}.
$$
So, if $\|f\|_{C(\mathbb T)}\leq 1$, then for every
$t\in\mathbb T$ we obtain
$$
|Q_N f(t)|=\bigg|m(t)\frac{1}{2\pi}\int_\mathbb T D_N(t-\theta)f(\theta)d\theta-
\frac{1}{2\pi}\int_\mathbb T D_N(t-\theta)m(\theta)f(\theta)d\theta\bigg|=
$$
$$
=\frac{1}{2\pi}\bigg|\int_\mathbb T D_N(t-\theta)(m(t)-m(\theta))f(\theta)d\theta\bigg|\leq
\frac{1}{2\pi}\int_{|t-\theta|\leq\pi} |D_N(t-\theta)||m(t)-m(\theta)|d\theta\leq c.
$$
Hence,
$$
\|Q_N\|_{C(\mathbb T)\rightarrow C(\mathbb T)}\leq c.
\eqno(3)
$$
Thus condition~(i) holds.

It is obvious that the product of two continuous functions
satisfying the uniform Dini condition (see (2)) satisfies the
uniform Dini condition as well and therefore belongs to
$U(\mathbb T)$. So, for all $n\in\mathbb Z$ we have $m e_n\in
U(\mathbb T)$, whence
$$
\|Q_N e_n\|_{C(\mathbb T)}=\|mS_N(e_n)-S_N(me_n)\|_{C(\mathbb
T)}\rightarrow 0 \quad \textrm{as} \,\, N\rightarrow\infty.
$$
Thus, (ii) holds.

It remains to note that from (i) and (ii) it follows that
$\|Q_N f\|_{C(\mathbb T)}\rightarrow 0$ for every function
$f\in C(\mathbb T)$, whence for every $f\in U(\mathbb T)$ we
obtain
$$
\|mf-S_N(mf)\|_{C(\mathbb T)}=\|m(f-S_N(f))+Q_N(f)\|_{C(\mathbb T)}\leq
$$
$$
\leq\|m\|_{C(\mathbb T)}\|f-S_N(f)\|_{C(\mathbb T)}+\|Q_N(f)\|_{C(\mathbb T)}\rightarrow 0.
$$
The bound of the norm $\|m\|_{MU(\mathbb T)}$ is obvious from
(3) and the identity $S_N(mf)=mS_N(f)-Q_N(f)$. This completes
the proof of Theorem 1 and, thereby, Theorems 2 and 3.

\quad

\begin{center}
\textbf{3. Two results of Salem type and their corollaries}
\end{center}

Recall Salem's result mentioned in the introduction: there
exist functions $f\in A(\mathbb T)$ and $g\in U(\mathbb T)$
such that $fg\notin U(\mathbb T)$. We note that the proof of
this result given in [1, Ch.~1, Sec.~6] yields in addition that
$\widehat{g}(k)=o(1/|k|)$. Modifying this proof, we obtain two
theorems of a similar type. As a consequence we shall see that
Theorem 3 is sharp.

Given a positive sequence $\gamma=\{\gamma(n), n=0, 1, 2,
\ldots\}$ consider the space $A_\gamma(\mathbb T)$ of all
functions $f\in C(\mathbb T)$ satisfying
$$
\|f\|_{A_\gamma(\mathbb T)}=\sum_{k\in\mathbb Z}|\widehat{f}(k)|\gamma(|k|)<\infty.
$$
We shall deal with sequences $\gamma$ bounded away from zero,
i.e., such that $\inf_n\gamma(n)>0$; this is why we assume the
continuity of functions in $A_\gamma(\mathbb T)$. If
$\gamma(n)\equiv 1$, we obtain the Wiener algebra $A(\mathbb
T)$. Clearly, $A_\gamma(\mathbb T)$ is a Banach space (provided
that $\gamma$ is bounded away from zero).

Let $V(\mathbb T)$ be the space of all functions of bounded
variation on $\mathbb T$. We denote the variation of $f$ on
$\mathbb T$ by $\|f\|_{V(\mathbb T)}$.

Consider also the space $W_2^{1/2}(\mathbb T)$ (the Sobolev
space) of all integrable functions $f$ on $\mathbb T$ with
$$
\|f\|_{W_2^{1/2}(\mathbb T)}=
\bigg(\sum_{k\in\mathbb Z}|\widehat{f}(k)|^2|k|\bigg)^{1/2}<\infty.
$$

Recall that $V\cap C(\mathbb T)\subseteq U(\mathbb T)$ and
$W_2^{1/2}\cap C(\mathbb T)\subseteq U(\mathbb T)$. These
inclusions are well known.\footnote{The former was obtained by
Jordan, see [4, Ch.~I, Sec.~39]. The latter can be proved as
follows. Assuming that $g\in W_2^{1/2}(\mathbb T)$, we have
$\sum_{|k|\leq N}|k||\widehat{g}(k)|=o(N)$, which for $g\in
C(\mathbb T)$ implies $g\in U(\mathbb T)$, see [4, Ch.~I,
Sec.~64].} It is also clear that if $\gamma$ is bounded away
from zero, then $A_\gamma(\mathbb T)\subseteq A(\mathbb
T)\subseteq U(\mathbb T)$.

\quad

\textsc{Theorem 4.} \emph{If
$$
\underset{n\rightarrow\infty}{\underline{\lim}}\,\frac{\gamma(n)}{\log n}=0,
$$
then there exist two (real) functions $f$ and $g$ such that
$f\in A_\gamma(\mathbb T)$ and $g\in V\cap C(\mathbb T)$ but
$fg\notin U(\mathbb T)$.}

\quad

\textsc{Theorem 5.} \emph{If
$$
\underset{n\rightarrow\infty}{\underline{\lim}}\,\frac{\gamma(n)}{(\log n)^{1/2}}=0,
$$
then there exist two (real) functions $f$ and $g$ such that
$f\in A_\gamma(\mathbb T)$ and $g\in W_2^{1/2}\cap V\cap
C(\mathbb T)$ but $fg\notin U(\mathbb T)$.}

\quad

Theorem 4 shows, in particular, that the logarithmic weight in
Theorem 3 cannot be replaced by a weight of slower growth.
Combining Theorems 3 and 4 we obtain the following corollary.

\quad

\textsc{Corollary 1.} \emph{The inclusion $A_\gamma(\mathbb
T)\subseteq MU(\mathbb T)$ holds if and only if
$$
\underset{n\rightarrow\infty}{\underline{\lim}}\,\frac{\gamma(n)}{\log n}>0.
$$}

\quad

The next corollary is an obvious consequence of Theorem 5.

\quad

\textsc{Corollary 2.} \emph{$W_2^{1/2}\cap V\cap C(\mathbb
T)\not\subseteq MU(\mathbb T)$.}

\quad

To prove Theorems 4 and 5 we need some lemmas.

It is convenient to consider the general case of linear normed
spaces $S$ embedded in $C(\mathbb T)$, i.e., such that $S$ is
contained in $C(\mathbb T)$ as a linear subspace and
$\|f\|_{C(\mathbb T)}\leq \mathrm{const} \|f\|_{S}$ for all
$f\in S$.

If $X$ and $Y$ are two sets in $C(\mathbb T)$ we let
$$
XY=\{xy :
x\in X, \,y\in Y\}.
$$

The following general lemma on multiplication is possibly known
in one form or another. We provide a short proof.

\quad

\textsc{Lemma 1.} \emph{Let $\mathbb X$, $\mathbb Y$, and
$\mathbb Z$ be Banach spaces embedded in $C(\mathbb T)$.
Suppose that $\mathbb X\mathbb Y\subseteq\mathbb Z$. Then there
exists a $c>0$ such that $\|xy\|_{\mathbb Z}\leq c
\|x\|_\mathbb X \|y\|_\mathbb Y$ for all $x\in\mathbb X$ and
$y\in\mathbb Y$.}

\quad

\textsc{Proof.} Given Banach spaces $E$ and $F$, let $B(E, F)$
denote the space of bounded operators from $E$ to $F$. For
$y\in\mathbb Y$, let $M_y \colon \mathbb X\rightarrow\mathbb Z$
be the operator that takes each element $x\in\mathbb X$ to the
product $xy$. Applying the closed graph theorem, we obtain
$M_y\in B(\mathbb X, \mathbb Z)$. Consider then the operator $Q
\colon \mathbb Y\rightarrow B(\mathbb X, \mathbb Z)$ that takes
each element $y$ of $\mathbb Y$ to the operator $M_y$. Applying
the closed graph theorem, we obtain $Q\in B(\mathbb Y,
B(\mathbb X, \mathbb Z))$. Setting $c=\|Q\|_{B(\mathbb Y,
B(\mathbb X, \mathbb Z))}$, we see that $ \|xy\|_{\mathbb
Z}=\|M_y x\|_{\mathbb Z}\leq \|M_y\|_{B(\mathbb X, \mathbb
Z)}\|x\|_{\mathbb X}= \|Qy\|_{B(\mathbb X, \mathbb
Z)}\|x\|_{\mathbb X}\leq c \|y\|_{\mathbb Y}\|x\|_{\mathbb X},
$. The lemma is proved.

\quad

For $n=1, 2, \ldots$, let $g_n$ be the functions defined by
$$
g_n(t)=\sum_{k=1}^n \bigg(1-\frac{k}{n}\bigg)\frac{1}{k}\sin kt.
\eqno(4)
$$

As above, by $e_n $ we denote the function
$e_n(t)=e^{int}$.

\quad

\textsc{Lemma 2.} \emph{The estimate $\|e_n
g_n\|_{U(\mathbb T)}\gtrsim\log n$ holds.}

\quad

\textsc{Proof.} We have
$$
g_n(t)=\sum_{1\leq |k|\leq n} \bigg(1-\frac{|k|}{n}\bigg)\frac{1}{2ik}e^{ikt},
$$
whence
$$
e_n(t) g_n(t)=\sum_{1\leq |k|\leq n} \bigg(1-\frac{|k|}{n}\bigg)\frac{1}{2ik}e^{i(k+n)t}.
$$
Therefore,
$$
S_n(e_n g_n)(t)=\sum_{1\leq |k|\leq n, \,\, |k+n|\leq n}\bigg(1-\frac{|k|}{n}\bigg)\frac{1}{2ik}e^{i(k+n)t}.
$$
This implies
$$
S_n(e_n g_n)(0)=\sum_{-n\leq k\leq -1}\bigg(1-\frac{|k|}{n}\bigg)\frac{1}{2ik}=
\sum_{k=1}^n\bigg(1-\frac{k}{n}\bigg)\frac{-1}{2ik}.
$$
So,
$$
\|e_n g_n\|_{U(\mathbb T)}\geq|S_n(e_n g_n)(0)|=\sum_{k=1}^n\bigg(1-\frac{k}{n}\bigg)\frac{1}{2k}=
\sum_{k=1}^n\frac{1}{2k}+O(1)\simeq\log n.
$$
The lemma is proved.

\quad

\textsc{Lemma 3.} \emph{The estimates $\|g_n\|_{V(\mathbb
T)}\leq 2\pi$ and $\|g_n\|_{C(\mathbb T)}\leq 2\pi$ hold.}

\quad

\textsc{Proof.} For the derivative $g_n'$ we have (see (4))
$$
g_n'(t)=\sum_{k=1}^n \bigg(1-\frac{k}{n}\bigg)\cos kt=\frac{1}{2}F_n(t)-\frac{1}{2},
$$
where $F_n$ is the Fej\'er kernel. As is known, $F_n(t)\geq 0$
for all $t$, so,
$$
|g_n'(t)|\leq \bigg|g_n'(t)+\frac{1}{2}\bigg|+\frac{1}{2}=
g_n'(t)+1,
$$
whence
$$
\|g_n\|_{V(\mathbb T)}=\int_0^{2\pi}|g_n'(t)| dt\leq 2\pi,
$$
which, in turn, implies $\|g_n\|_{C(\mathbb T)}\leq 2\pi$,
because $g_n(0)=0$. The lemma is proved.

\quad

\textsc{Proof of Theorems 4 and 5.} We set
$$
\|f\|_{V\cap C(\mathbb T)}=\|f\|_{V(\mathbb T)}+\|f\|_{C(\mathbb T)},
$$
$$
\|f\|_{W_2^{1/2}\cap V\cap C(\mathbb T)}=\|f\|_{W_2^{1/2}(\mathbb T)}+
\|f\|_{V(\mathbb T)}+\|f\|_{C(\mathbb T)}.
$$
The spaces $V\cap C(\mathbb T)$ and $W_2^{1/2}\cap V\cap
C(\mathbb T)$ endowed with these norms are Banach spaces (we
leave a simple verification to the reader).

Thus, each of the four spaces $U(\mathbb T)$,
$A_\gamma(\mathbb T)$, $V\cap C(\mathbb T)$, and $W_2^{1/2}\cap
V\cap C(\mathbb T)$ is a Banach space
embedded in $C(\mathbb T)$.

We apply Lemma 1. Assuming that the conclusion of Theorem 4 is
false, we obtain $\|fg\|_{U}\leq c\|f\|_{A_\gamma}\|g\|_{V\cap
C}$. In particular,
$$
\|e_ng_n\|_{U}\lesssim\|e_n\|_{A_\gamma}\|g_n\|_{V\cap C}.
\eqno(5)
$$
Clearly, $\|e_n\|_{A_\gamma}=\gamma(n)$. Taking Lemmas 2 and 3
into account and using (5), we obtain $\log
n\lesssim\gamma(n)$, which contradicts the condition of Theorem
4.

Assuming that the conclusion of Theorem 5 is false, we have
$$
\|e_ng_n\|_{U}\lesssim\|e_n\|_{A_\gamma}\|g_n\|_{W_2^{1/2}\cap V\cap C}.
\eqno(6)
$$
Obviously, $\|g_n\|_{W_2^{1/2}}\simeq(\log n)^{1/2}$, whence,
taking Lemma 3 into account, we obtain $\|g_n\|_{W_2^{1/2}\cap
V\cap C}\simeq(\log n)^{1/2}$. Thus, relation (6) and Lemma 2
yield $\log n\lesssim\gamma(n)(\log n)^{1/2}$, which
contradicts the condition of Theorem 5.

\quad

\textbf{Remarks.} 1. The following result, directly related to
the topic of this paper, was obtained by Olevskii [5]: The
algebra generated by $U(\mathbb T)$ coincides with $C(\mathbb
T)$; moreover, each function $f\in C(\mathbb T)$ can be
represented in the form $f=\varphi_1\varphi_2+\varphi_3$, where
$\varphi_j\in U(\mathbb T)$, $j=1, 2, 3$. Apparently, it is
unknown whether $\varphi_3$ can be set to zero. It would be
interesting to find out if one can chose the functions
$\varphi_j$ so that they belong to more narrow classes rather
then $U(\mathbb T)$; for instance, can we chose $\varphi_1\in
A(\mathbb T)$?

2. Theorem 3 immediately implies $\|e_n\|_{MU}\lesssim\log
|n|$. At the same time it is easy to see that the functions
$g_n$ defined by (4) satisfy $\|g_n\|_{U}=O(1)$. So, using
Lemma 2, we see that $\log|n|\lesssim\|e_ng_n\|_{U}
\leq\|e_n\|_{MU}\|g_n\|_{U}\lesssim\|e_n\|_{MU}$. Thus,
$\|e_n\|_{MU(\mathbb T)}\simeq\log |n|$.

3. It is natural to consider the nonsymmetric analogue of the
space $U(\mathbb T)$, namely, the space $U^{asym}(\mathbb T)$
defined in the similar way as $U(\mathbb T)$ with the only
difference that instead of the partial sums $S_N(f)(t)$ one
uses the partial sums $S_{N, M}(f)(t)=\sum_{-N\leq k\leq
M}\widehat{f}(k)e^{ikt}$. The space $MU^{asym}(\mathbb T)$ of
multipliers is defined in a natural way. Note that the spaces
of multipliers in the symmetric and nonsymmetric cases are
different. It is easy to verify that
$\|e_n\|_{MU^{asym}(\mathbb T)}=O(1)$ and, therefore,
$A(\mathbb T)\subseteq MU^{asym}(\mathbb T)$. It is worth
mentioning that this embedding has a counterpart for functions
analytic in the disk $D=\{z\in \mathbb C : |z|<1\}$ of the
complex plane $\mathbb C$. Let $U^+(D)$ be the class of
functions analytic in $D$ whose Taylor series converges
uniformly in $D$, and let $A^+(D)$ be the class of functions
analytic in $D$ whose sequence of Taylor coefficients belongs
to $l^1$. If $f\in A^+(D)$ and $g\in U^+(D)$, then $fg\in
U^+(D)$.

4. The conditions in Theorems 1--3 of this paper are resemblant
to the conditions that appear in the work by Vinogradov,
Goluzina, and Khavin [6, Theorem 3] on multipliers of the space
$K$ of Cauchy--Stieltjes-type integrals. The resemblance of
these conditions is a reflection of a certain likeness between
the character of singularity of the Dirichlet and Cauchy
kernels.

\quad

\begin{center}
\textbf{References}
\end{center}

\flushleft
\begin{enumerate}

\item J.-P. Kahane, \emph{S\'erie de Fourier absolument
    convergentes}, Springer-Verlag, Berlin--Heidelberg--New
    York, 1970.

\item  A. N. Kolmogorov and S. V. Fomin, \emph{Elements of the
    Theory of Functions and Functional Analysis}, Nauka,
    Moskov 1989; English transl.: Dover Publications, Inc.,
    Mineola, New York, 1999.

\item  V. V. Lebedev, ``On Uniform Convergence of Fourier
    Series,'' Mat. Zametki, \textbf{91}:6 (2012), 946--949;
    English transl.: Math. Notes, \textbf{91}:6 (2012),
    889--892.

\item  N. K. Bari, \emph{A treatise on trigonometric series},
    Fizmatgiz, Moscow 1961; English transl.: Pergamon, New
    York 1964.

\item A. M. Olevskii, ``On the algebra of functions,
    generated by uniformly convergent Fourier series,'' DAN
    SSSR, \textbf{297}:4 (1987), 798--800; English transl.:
    Soviet Math. Dokl. \textbf{36}:3 (1988), 542--544.

\item S. A. Vinogradov, M. G. Goluzina, and V. P. Khavin,
    ``Multipliers and divisors of Cauchy--Stiltjes
    integrals,'' Zap. Nauchn. Sem. LOMI, 19 (1970), 55--78
    (in Russian).

\end{enumerate}

\quad

\quad

\qquad National Research University Higher School of Economics

\qquad e-mail address: \emph{lebedevhome@gmail.com}

\end{document}